\documentclass[11pt]{amsart}
\usepackage{amsmath}
\usepackage{amssymb}
\usepackage{amsfonts}
\usepackage{amsthm}
\usepackage{enumerate}
\usepackage{enumitem}
\usepackage[mathscr]{eucal}
\usepackage{verbatim}
\usepackage{amsthm}
\usepackage{amscd}
\usepackage[mathscr]{eucal}
\usepackage{appendix}
\usepackage{tikz}

\numberwithin{equation}{section}

\newtheorem{innercustomgeneric}{\customgenericname}
\providecommand{\customgenericname}{}

\newcommand{\newcustomtheorem}[2]{\newenvironment{#1}[1]
  {\renewcommand\customgenericname{#2}
   \renewcommand\theinnercustomgeneric{##1}\innercustomgeneric}{\endinnercustomgeneric}}

\newcustomtheorem{customthm}{Theorem}

\newcommand{\newcustomlemma}[2]{\newenvironment{#1}[1]
  {\renewcommand\customgenericname{#2}
   \renewcommand\theinnercustomgeneric{##1} \innercustomgeneric}{\endinnercustomgeneric}}

\newcustomlemma{customlemma}{Lemma}

\newcustomlemma{customproposition}{Proposition}

\newcustomlemma{customclaim}{Claim}

\theoremstyle{plain}
\newtheorem{theorem}{Theorem}
\newtheorem{corollary}[theorem]{Corollary}
\newtheorem{lemma}[theorem]{Lemma}

\newtheorem*{theorem*}{Theorem}
\newtheorem*{lemma*}{Lemma}
\newtheorem*{proposition*}{Proposition}
\newtheorem*{corollary*}{Corollary}
\newtheorem*{remark*}{Remark} 
\newtheorem*{remarks*}{Remarks}
\newtheorem*{conj*}{Conjecture}

\newcommand{\bbz}{\mathbb{Z}}

\newcommand{\bbr}{\mathbb{R}}
\newcommand{\bbrn}{\mathbb R^n}

\newcommand{\bbn}{\mathbb{N}}

\newcommand{\xxxi}{\vec{\boldsymbol{\xi}\,}}

\newcommand{\yyy}{\vec{\boldsymbol{y}}}
\newcommand{\www}{\vec{\boldsymbol{w}}}

\newcommand{\ppp}{\vec{\boldsymbol{p}}}

\newcommand{\uuu}{\vec{\boldsymbol{u}}}

\def\000{\vec{\boldsymbol{0}}}

\newcommand{\q}{\quad}

\DeclareFontFamily{U}{mathx}{\hyphenchar\font45}
\DeclareFontShape{U}{mathx}{m}{n}{
	<5> <6> <7> <8> <9> <10>
	<10.95> <12> <14.4> <17.28> <20.74> <24.88>
	mathx10
}{}

\def\wh{\widehat}

\newcommand{\supp}{\mathrm{supp}}

\makeatletter
\@namedef{subjclassname@2020}{\textup{2020} Mathematics Subject Classification}
\makeatother

\allowdisplaybreaks

\makeatletter
\@namedef{subjclassname@2020}{\textup{2020} Mathematics Subject Classification}
\makeatother

\makeindex         

\begin{document}

\author{Bae Jun Park}
\address{B. Park, Department of Mathematics, Sungkyunkwan University, Suwon 16419, Republic of Korea}
\email{bpark43@skku.edu}

\author{Naohito Tomita}
\address{N. Tomita, Graduate School of Science, The University of Osaka, Toyonaka, Osaka 560-0043, Japan}
\email{tomita@math.sci.osaka-u.ac.jp}

\thanks{B. Park is supported by the National Research Foundation (Republic of
Korea) grant RS-2025-20512969 and by the Open KIAS Center at Korea Institute for Advanced Study.   N. Tomita is supported by JSPS KAKENHI Grant Number 25K07086.}

 \title[Necessary conditions for weighted estimates of Multilinear operators]{Necessary conditions for weighted estimates of Multilinear Multipliers and Pseudo-Differential Operators} 
\subjclass[2020]{Primary 42B20, 42B25, 47H60}
\keywords{Weighted estimates, Muckenhoupt weights, Multiple weight spaces, Multilinear multipliers, Multilinear pseudo-differential operators }

\begin{abstract}
We study optimal multiple weight assumptions in the weighted theory of multilinear Fourier multipliers and multilinear pseudo-differential operators.
For multilinear Fourier multipliers, we revisit the weighted H\"ormander-type theorem of Li and Sun \cite{Li_Su2015}, as a multilinear version of Kurtz and Wheeden \cite{Ku_Wh1979}, and show that their multiple weight condition is sharp. This provides the sharp necessary condition in the multilinear setting and simultaneously improves the classical linear necessity established by Kurtz and Wheeden \cite{Ku_Wh1981}.
In the pseudo-differential setting, we consider recent weighted estimates of the authors \cite{Park_To2024, Park_To2025} for symbols in the multilinear H\"ormander class and prove that their multiple weight hypothesis is also best possible.
As a corollary, we can obtain the optimality of sharp maximal function estimates for multilinear pseudo-differential operators in \cite{Park_To2024, Park_To2025} which originated in Chanillo and Torchinsky \cite{Ch_To1985}.
\end{abstract}

\maketitle

%\tableofcontents
\section{Introduction}\label{introsec}

Let $w \colon \mathbb{R}^n \to [0,\infty)$ be a locally integrable function. 
For $1<p<\infty$, we say that $w$ belongs to the Muckenhoupt class $A_p(\bbrn)$ if
\begin{equation*}
    [w]_{A_p}
    := \sup_{Q}
    \left( \frac{1}{|Q|} \int_Q w(x)\,dx \right)^{\frac{1}{p}}
    \left( \frac{1}{|Q|} \int_Q \big( w(x) \big)^{1-p'}\,dx \right)^{\frac{1}{p'}}
    < \infty,
\end{equation*}
where  $p' = \frac{p}{p-1}$ denotes the H\"older conjugate of $p$, and
the supremum is taken over all cubes $Q \subset \mathbb{R}^n$ whose sides
are parallel to the coordinate axes.
The quantity $[w]_{A_p}$ is called the $A_p$ constant of the weight $w$.
For $p=1$, the class $A_1(\bbrn)$ consists of all weights satisfying
\begin{equation*}
    \mathcal{M}w(x) \le C\, w(x)
    \qquad \text{for a.e. } x\in\mathbb{R}^n,
\end{equation*}
where $\mathcal{M}$ denotes the Hardy--Littlewood maximal operator taken over 
axis-parallel cubes, and the smallest such constant $C$ is denoted by $[w]_{A_1}$.
Equivalently, $w\in A_1$ if
\begin{equation*}
    \frac{1}{|Q|} \int_Q w(y)\,dy
    \;\le\;
    [w]_{A_1} \inf_{x\in Q} w(x)
    \qquad \text{for all cubes } Q\subset\mathbb{R}^n.
\end{equation*}
Then it turns out that 
 \begin{equation}\label{apweightinclusion}
 A_p(\bbrn)\subset A_q(\bbrn) \q  \text{ if}~ 1\le p\le q< \infty.
 \end{equation}
Given a weight $w$,
the weighted Lebesgue space $L^p(w)$, $1<p<\infty$, consists of all measureble functions $f$ on $\bbrn$ satisfying
$$\|f\|_{L^p(w)}:=\left(\int_{\bbrn}|f(x)|^p w(x)\, dx \right)^{\frac{1}{p}}<\infty . $$
A classical result of Muckenhoupt \cite{Mu1972} asserts the fundamental equivalence
\begin{equation}\label{chaweight-intro}
\|\mathcal{M}f\|_{L^p(w)}\lesssim \|f\|_{L^p(w)}
\quad\Longleftrightarrow\quad
w\in A_p(\bbrn),
\qquad(1<p<\infty),
\end{equation}
which forms the foundation of the modern weighted theory for Calder\'on-Zygmund operators.
Indeed, it is well-known that for all $1<p<\infty$,
\begin{equation}\label{czweightest}
\text{ any Calder\'on-Zygmund operators $T$ are bounded on $L^p(w)$ if $w\in A_p(\bbrn)$.}
\end{equation}

\medskip

As the field progressed toward a deeper understanding of multilinear operators, it became necessary to develop a suitable analogue of the $A_p(\bbrn)$ theory for tuples of weights.  
This multilinear extension was introduced by Lerner, Ombrosi, P\'erez, Torres, and Trujillo-Gonz\'alez \cite{Le_Om_Pe_To_Tr2009}, who defined the multiple weight class $A_{\ppp}(\bbrn)=A_{(p_1,\dots,p_l)}(\bbrn)$ associated with exponents $1\le p_1,\dots,p_l<\infty$. 
An $l$-tuple of weights $\www=(w_1,\dots,w_l)$, defined on $\bbrn$, belongs to $A_{\ppp}(\bbrn)$ if
\begin{equation}\label{multiplewedef}
[\www]_{A_{\ppp}}:=\sup_Q
\Bigg[
\left(\frac{1}{|Q|}\int_Q v_{\www}(x)\,dx\right)^{\frac{1}{p}}
\prod_{j=1}^l
\left(\frac{1}{|Q|}\int_Q \big(w_j(x)\big)^{\,1-p_j'}\,dx\right)^{\frac{1}{p_j'}}
\Bigg]<\infty,
\end{equation}
where  the exponent $p$ is determined by $\frac1p=\frac1{p_1}+\cdots+\frac1{p_l}$ and 
$$v_{\www}(x)=\prod_{j=1}^l w_j(x)^{\frac{p}{p_j}}.$$ 
When $p_j=1$, $\left(\frac{1}{|Q|}\int_Q \big(w_j(x)\big)^{\,1-p_j'}\,dx\right)^{\frac{1}{p_j}}$ is understood as $\big(\inf_Q w_j\big)^{-1}$.

One of the most striking achievements of \cite{Le_Om_Pe_To_Tr2009} is the discovery that $A_{\ppp}(\bbrn)$ can be characterized through a maximal inequality fully analogous to the linear characterization \eqref{chaweight-intro}.  
For the multi-sublinear Hardy-Littlewood maximal operator
\[
\mathbf{M}(f_1,\dots,f_l)(x)
  :=\sup_{Q\ni x}
      \Bigg(
        \frac{1}{|Q|^l}
        \int_{Q^l}
          \prod_{j=1}^l |f_j(u_j)|\, d\uuu
      \Bigg)
\]
where $Q^l:=Q\times \cdots \times Q$, $d\uuu:=du_1\cdots du_l$, and
 the supremum is taken over all cubes $Q$ in $\bbrn$ containing $x$,
it was shown in \cite{Le_Om_Pe_To_Tr2009} that
\begin{equation}\label{Le_Om_Pe_To_Tr2009_result}
\big\Vert \mathbf{M}(f_1,\dots,f_l)\big\Vert_{L^p(v_{\www})}
   \lesssim
   \prod_{j=1}^l \|f_j\|_{L^{p_j}(w_j)}
   \quad\Longleftrightarrow\quad
   \www\in A_{\ppp}(\bbrn)
\end{equation}
with usual modifications if some $p_j$ are one.
This result identifies the class $A_{\ppp}(\mathbb{R}^n)$ as the natural framework for the weighted theory of multilinear harmonic analysis.
However, the internal structure of $A_{\ppp}(\mathbb{R}^n)$ is considerably more intricate than in the linear case.
For instance, these classes do not satisfy the monotonicity property \eqref{apweightinclusion} known for the classical $A_p$ weights, and the interaction among the component weights is rather subtle; see \cite[Remark~7.3]{Le_Om_Pe_To_Tr2009}.
A key decomposition theorem in \cite{Le_Om_Pe_To_Tr2009} shows that any multiple weight $\www\in A_{\ppp}(\mathbb{R}^n)$ can be characterized entirely in terms of scalar $A_q$ conditions: each component weight $w_j$ satisfies a suitable linear Muckenhoupt condition, and the combined weight $v_{\www}$ lies in another scalar class (see Lemma~\ref{multiindiweight}).
This decomposition clarifies the geometry of $A_{\ppp}(\mathbb{R}^n)$ and will be crucial for our purposes, particularly in constructing counterexamples that establish the sharpness of our main results.

\bigskip

%%%%%%%%%%%%%%%%%%%%%%%%%%%%%%%%%%%%%%%%%%%%%%%%%%%%
% 1.  MULTILINEAR MULTIPLIERS
%%%%%%%%%%%%%%%%%%%%%%%%%%%%%%%%%%%%%%%%%%%%%%%%%%%%
\subsection*{Multilinear Fourier multipliers}
In the linear setting, 
let $\sigma$ be a bounded function on $\bbrn$ and we define the corresponding multiplier operator $T_{\sigma}$ as
\begin{equation}\label{limulopdef}
T_\sigma f(x)
  = \int_{\mathbb{R}^n}  \sigma(\xi)\widehat{f}(\xi) e^{2\pi i \langle x,\xi\rangle}
   \, d\xi
\end{equation}
for all Schwartz functions $f$ in $\bbrn$.  
The classical theory, dating back to Mikhlin \cite{Mik1956} and H\"ormander \cite{Ho1960}, shows that if
\begin{equation*}
\sup_{k\in\mathbb{Z}}
\|\sigma(2^k\cdot)\widehat{\psi}\|_{L^2_s(\bbrn)}<\infty,
\qquad s>\frac{n}{2},
\end{equation*}
then $T_\sigma$ is bounded on $L^p(\bbrn)$ for every $1<p<\infty$.  
Here, $\psi$ is a Schwartz function on $\bbrn$ whose Fourier transform is supported in the annulus $\{\xi\in\bbrn: 2^{-1}\le |\xi|\le 2\}$ and satisfies $\sum_{j\in\bbz}\wh{\psi}(2^j\xi)=1$ for $\xi\not=0$,
and the norm of the Sobolev space $L^2_s(\bbrn)$ is defined by
\[
\|f\|_{L^2_s(\bbrn)}
=\left(\int_{\bbrn}(1+4\pi^2|\xi|^2)^s|\widehat{f}(\xi)|^2\, d\xi \right)^{\frac{1}{2}}.
\]
The regularity condition $s>\frac{n}{2}$ is known to be optimal for the $L^p$ boundedness to hold.  
Weighted extensions of this theorem were developed by Kurtz and Wheeden \cite{Ku_Wh1979}, establishing that the same Sobolev-type regularity condition suffices for boundedness on weighted spaces $L^p(w)$ with $w\in A_{\frac{ps}{n}}(\bbrn)$ for $p>\frac{n}{s}$.
   \begin{customthm}{A}\cite{Ku_Wh1979}\label{linearmultiplier}
Suppose that
\begin{equation}\label{hormultiass}
\sup_{k\in\bbz}\big\Vert \sigma(2^k\cdot)\widehat{\psi}\big\Vert_{L^2_s(\bbrn)}<\infty, \quad \frac{n}{2}<s\le n.
\end{equation}
For any $\frac{n}{s}<p<\infty$ and $w\in A_{\frac{ps}{n}}(\bbrn)$,
we have
\begin{equation}\label{tsiglpwboun}
\big\Vert T_{\sigma}\big\Vert_{L^p(w)\to L^p(w)}\lesssim_w \sup_{k\in\bbz}\big\Vert \sigma(2^k\cdot)\widehat{\psi}\big\Vert_{L^2_s(\bbrn)}.
\end{equation}
 \end{customthm}

In light of \eqref{czweightest}, one might speculate that, while keeping all other assumptions in Theorem $\ref{linearmultiplier}$ unchanged, the weaker requirement $w\in A_p(\bbrn)$ could still suffice to ensure the $L^p(w)$ boundedness of $T_{\sigma}$. However, Kurtz and Wheeden \cite{Ku_Wh1981} proved that this is not the case.
\begin{customthm}{B}\cite{Ku_Wh1981}\label{linearsharp}
Let  $s$ be an integer satisfying $\frac{n}{2}<s<n$. Then for any $\frac{n}{s}\le p\le (\frac{n}{s})'$, there exist $w\in A_{p}(\bbrn)$ and $\sigma$ satisfying \eqref{hormultiass} such that \eqref{tsiglpwboun} does not hold.
\end{customthm}   
Although this shows that $A_p$ weights are not strong enough to guarantee the weighted estimate, they are still much weaker than $ A_{\frac{ps}{n}}$ weights required in Theorem \ref{linearmultiplier}, where the exponent $\frac{ps}{n}$ is strictly smaller than $p$.  Thus, the question of identifying the optimal weight condition remains open.

In this paper, we will completely  address this issue by showing that the condition $w\in A_{\frac{ps}{n}}(\bbrn)$ is indeed optimal for the weighted $L^p$ boundedness. In fact, we will establish this optimality in a broader framework by proving it in the multilinear setting. (see Theorem \ref{sharpnessmultilinear} below)

\medskip

The multilinear generalization of multiplier theory was initiated by Coifman and Meyer \cite{Co_Me1978}, who introduced operators of the form
\[
T_\sigma(f_1,\dots,f_l)(x)
 = \int_{(\mathbb{R}^n)^l}
     \sigma(\xi_1,\dots,\xi_l)
     \Big( \prod_{j=1}^l \widehat{f_j}(\xi_j)\Big)
     e^{2\pi i \langle x,\xi_1+\cdots+\xi_l\rangle}\, d\xxxi
\]
where $d\xxxi:=d\xi_1\cdots d\xi_l$.
Their work generalized the classical result of Mihklin \cite{Mik1956}, establishing $L^{p_1}\times\cdots\times L^{p_l}\to L^p$ boundedness for $1<p_1,\dots,p_l<\infty$ and $1\le p<\infty$ with $\frac{1}{p_1}+\cdots+\frac{1}{p_l}=\frac{1}{p}$, provided that the symbol $\sigma$ satisfies a sufficiently large regularity condition. The boundedness result was later extended to all parameters $\frac{1}{l}<p<\infty$ by Kenig and Stein \cite{Ke-St1999} and Grafakos and Torres \cite{Gr_To2002}. Afterwards, Tomita \cite{To2010} provided a sharp estimate for $T_{\sigma}$ as a multilinear counterpart of H\"ormander's result with $1<p<\infty$ and this was later extended by Grafakos and Si \cite{Gr_Si2012} to the range $\frac{1}{l}<p\le 1$. Let $\Psi$ be a Schwartz function on $\bbr^{nl}$, which is a multilinear counterpart of $\psi$, satisfying
 $\supp{\wh{\Psi}}\subset \{\xxxi\in (\bbrn)^l: 2^{-1}\le |\xxxi|\le 2\}$ and  $\sum_{j\in\bbz}\wh{\Psi}(2^j\xxxi)=1$ for $\xxxi\not=\vec{0}$.

\begin{customthm}{C} [\cite{Gr_Si2012, To2010}]\label{unweighted-multilinear-intro}
Let $1<p_1,\dots,p_l<\infty$ with $\frac{1}{p_1}+\dots+\frac{1}{p_l}= \frac{1}{p}$.  
Suppose that a multilinear symbol $\sigma$ satisfies
$$ \sup_{k\in\mathbb{Z}}
\big\Vert \sigma(2^k\vec{\cdot}\,)\widehat{\Psi}\big\Vert_{L^2_s(\bbr^{nl})}<\infty,
\qquad
s>\frac{nl}{2}.
$$
For any $\frac{nl}{s}<p_1,\dots,p_l<\infty$, we have
\begin{equation}\label{unweimulintest}
\big\Vert T_{\sigma}\big\Vert_{L^{p_1}(\mathbb{R}^n)\times\cdots\times L^{p_l}(\mathbb{R}^n)\to L^p(\bbrn)}\lesssim \sup_{k\in\mathbb{Z}}
\big\Vert \sigma(2^k\vec{\cdot}\,)\widehat{\Psi}\big\Vert_{L^2_s(\bbr^{nl})}.
\end{equation}
\end{customthm}
It is known that the estimate \eqref{unweimulintest} in Theorem~\ref{unweighted-multilinear-intro}
holds for a wider range of exponents \(p_1,\dots,p_l\),
while the present paper is concerned only with the case
\(\frac{nl}{s}<p_1,\dots,p_l<\infty\),
and the theorem is therefore stated in the above form.
For a precise description of the full admissible range of exponents,
we refer the reader to \cite{Gr_Si2012, To2010} for further details,
and to \cite{LHHLPPY} for a version covering the full range
\(0<p_1,\dots,p_l<\infty\).

\medskip

There exist various important results concerning weighted multilinear multiplier theorems, and among these we highlight two representative contributions that are most relevant to the present work.  
Fujita and Tomita \cite{Fu_To2012} obtained weighted estimates for the multilinear multiplier theorem stated in Theorem \ref{unweighted-multilinear-intro} by employing the classical Muckenhoupt $A_p(\bbrn)$ weights.  
Later, Li and Sun \cite{Li_Su2015} refined and improved this direction by establishing weighted estimates for the same multilinear multiplier results within the full multiple weight framework.
We also refer to the work of Bui and Duong \cite{Bu-Du2013} for related results in the bilinear setting.
\begin{customthm}{D}\cite{Li_Su2015}\label{multmulweight}
Let $1<p_1,\dots,p_l<\infty$ with $\frac{1}{p_1}+\dots+\frac{1}{p_l}= \frac{1}{p}$. 
Suppose that a multilinear symbol $\sigma$ satisfies
\[
\sup_{k\in\mathbb{Z}}
\|\sigma(2^k\vec{\cdot}\,)\widehat{\Psi}\|_{L^2_s(\bbr^{nl})}<\infty, \quad \frac{nl}{2}<s\le nl.
\]
For any $\frac{nl}{s}<p_1,\dots,p_l<\infty$ and $\www\in A_{(\frac{p_1s}{nl},\dots,\frac{p_l s}{nl})}(\bbrn)$,
we have
\begin{equation}\label{multiweibound}
\big\|T_\sigma\big\|_{L^{p_1}(w_1)\times \cdots\times L^{p_l}(w_l)\to L^p(v_{\www})}
\lesssim_{\www}
\sup_{k}\big\Vert \sigma(2^k\vec{\cdot}\,)\widehat{\Psi}\big\Vert_{L^2_s(\bbr^{nl})}.
\end{equation}
\end{customthm}

\medskip

Our first main theorem shows that the weighted estimate above cannot hold for larger weight classes $A_{(\frac{p_1}{q},\dots,\frac{p_l}{q})}(\bbrn)$ with $q< \frac{nl}{s}$.  
We remark that when $l\ge 2$, the classes $A_{\ppp}(\bbrn)$ are not generally increasing with the natural partial order on exponents $(p_1,\dots,p_l)$, unlike \eqref{apweightinclusion}, but they satisfy an increasing property under uniform scaling: for the scaled exponents $(rp_1,\dots,rp_l)$ with $r>1$, the classes $A_{(rp_1,\dots,rp_l)}$ become larger as $r$ increases. See \cite[Lemma 5]{Park2024_submitted} for more details.

\begin{theorem}\label{sharpnessmultilinear}
Let $l\in \bbn$, $\frac{nl}{2}<s\le nl$, $0<q<\infty$, $\max\{q,\frac{nl}{s}\} \le p_1,\dots,p_l<\infty$,
and $\frac{1}{p_1}+\cdots+\frac{1}{p_l}=\frac{1}{p}$.
Assume that \eqref{multiweibound} holds for all multilinear symbols $\sigma$ and for all $\vec{w}\in A_{(\frac{p_1}{q},\dots,\frac{p_l}{q})}(\bbrn)$.
Then it is necessary that
\[
q \ge \frac{nl}{s}.
\]
\end{theorem}

This result not only improves Theorem \ref{linearsharp} in the linear case $l=1$, but also provides its full multilinear analogue.

\bigskip

%%%%%%%%%%%%%%%%%%%%%%%%%%%%%%%%%%%%%%%%%%%%%%%%%%%%
% 2. MULTILINEAR PSEUDO-DIFFERENTIAL OPERATORS
%%%%%%%%%%%%%%%%%%%%%%%%%%%%%%%%%%%%%%%%%%%%%%%%%%%%

\subsection*{Multilinear pseudo-differential operators}

Given $m\in\mathbb{R}$ and $0\le\delta\le\rho\le 1$, the H\"ormander symbol class $S_{\rho,\delta}^m(\mathbb{R}^n)$ consists of all smooth functions $a(x,\xi)$ on $\mathbb{R}^n\times\mathbb{R}^n$ such that
\[
\big|\partial_x^\alpha\partial_\xi^\beta a(x,\xi) \big|
   \lesssim (1+|\xi|)^{m+\delta|\alpha|-\rho|\beta|}
\qquad x,\xi\in\bbrn
\]
for any multi-indices $\alpha,\beta \in\mathbb{N}_0^n$.
For $a \in S_{\rho,\delta}^m(\mathbb{R}^n)$, the associated (linear) pseudo-differential operator is defined, initially on the Schwartz class $\mathscr{S}(\mathbb{R}^n)$, by
\[
T_{[a]}f(x)
   = \int_{\mathbb{R}^n} a(x,\xi)\widehat{f}(\xi)
      e^{2\pi i \langle x,\xi\rangle}\,d\xi.
\]
 Denote by $\mathrm{Op}S_{\rho,\delta}^{m}(\bbrn)$ the class of pseudo-differential operators with symbols in $S_{\rho,\delta}^{m}(\bbrn)$.
Compared to \eqref{limulopdef}, when $a$ is independent of $x$, this reduces to a (linear) Fourier multiplier operator, so the class $\mathrm{Op}\,S_{\rho,\delta}^m(\mathbb{R}^n)$ can be viewed as a variable-coefficient generalization of the multiplier operators considered earlier.

The boundedness properties of such operators were developed systematically by H\"ormander \cite[Theorem 3.5]{Ho1967}, Calder\'on and Vaillancourt \cite{Ca_Va1972}, and Fefferman \cite{Fe1973}.  
In particular, for $0\le \delta\le \rho<1$ and $1<p<\infty$, the condition
\[
m \le -n(1-\rho)\Bigl|\frac12-\frac1p\Bigr|
\]
is sufficient for $L^p$-boundedness when the symbol belongs to $S_{\rho,\delta}^m(\mathbb{R}^n)$.  
This degree condition is often regarded as the linear prototype for many subsequent pseudo-differential estimates.

\medskip

Weighted estimates for linear pseudo-differential operators were initiated by Chanillo and Torchinsky \cite{Ch_To1985}, who proved if $a\in S_{\rho,\delta}^{-\frac{n}{2}(1-\rho)}(\bbrn)$ for $0\le \delta<\rho<1$ and $w\in A_{\frac{p}{2}}(\bbrn)$ for $2\le p<\infty$, then $T_{[a]}$ is bounded in $L^p(w)$.
Their proof relies essentially on a pointwise control of $T_{[a]}f$ by a sharp maximal function.  
This approach was later refined by Miyachi and Yabuta \cite{Mi_Ya1987}, who extended the underlying pointwise estimate to an $L^r$-sharp maximal function with $1<r\le 2$, thus providing a more flexible framework for weighted inequalities.  
Recently, the authors \cite{Park_To2024, Park_To2025} established a weighted estimate that unifies all the results discussed above. Since this result applies not only to the linear case but also to the general multilinear setting, we first describe some boundedness result for multilinear pseudo-differential operators and then state the result of \cite{Park_To2024, Park_To2025}.

\medskip

The multilinear H\"ormander symbol class $\mathbb{M}_lS_{\rho,\delta}^{m}(\bbrn)$ consists of all smooth functions $a$ on $(\bbrn)^{l+1}$ having the property that for all multi-indices $\alpha, \beta_1,\dots,\beta_l\in (\bbn_0)^n$ there exists a constant $C_{\alpha,\beta_1,\dots,\beta_l}>0$ such that
$$\big| \partial_{x}^{\alpha}\partial_{\xi_1}^{\beta_1}\cdots\partial_{\beta_l}^{\beta_l}a(x,\xi_1,\dots,\xi_l)\big|\le C_{\alpha,\beta_1,\dots,\beta_l}\big( 1+|\xi_1|+\cdots+|\xi_l|\big)^{m+\delta |\alpha|-\rho(|\beta_1|+\dots+|\beta_l|)},$$
 and let $T_{[a]}$ now denote the multilinear pseudo-differential operator associated with $a\in \mathbb{M}_lS_{\rho,\delta}^{m}(\bbrn)$, defined by
 $$T_{[a]}\big(f_1,\dots,f_l\big)(x):=\int_{(\bbrn)^l}a(x,\xi_1,\dots,\xi_l)\Big( \prod_{j=1}^{l}\wh{f_j}(\xi_j)\Big) \, e^{2\pi i\langle x,\xi_1+\dots+\xi_l\rangle}\; d\xxxi$$
for $f_1,\dots,f_l\in \mathscr{S}(\bbrn)$. 
In contrast with the $L^2$ boundedness of linear operators associated with $a\in S_{0,0}^{0}(\bbrn)$ in \cite{Ca_Va1972}, B\'enyi and Torres \cite{Be_To2004} showed that the bilinear operators associated with $a\in\mathbb{M}_2S_{0,0}^0(\mathbb{R}^n)$ need not be bounded on any reasonable product of Lebesgue spaces.  
This failure of boundedness has motivated an extensive study of the exact degree conditions needed for bilinear and multilinear pseudo-differential operators.

A symbolic calculus for bilinear operators was developed by B\'enyi, Maldonado, Naibo, and Torres \cite{Be_Ma_Na_To2010}, and a series of works by Michalowski, Rule, and Staubach \cite{Mi_Ru_St2014}, and  B\'enyi, Bernicot, Maldonado, Naibo, and Torres \cite{Be_Be_Ma_Na_To2013} established the  $L^{p_1}\times L^{p_2}\to L^p$ bounds for $1\le p,p_1,p_2\le \infty$ in the subcritical regime $m<m_\rho(p_1,p_2)$.
More recently, Miyachi and Tomita \cite{Mi_To2013, Mi_To2019, Mi_To2020} identified the optimal critical order $m=m_\rho(p_1,p_2)$ by obtaining a sharp bilinear counterpart of the $L^p$ boundedness of Fefferman \cite{Fe1973}.  
For general multilinear operators, Kato, Miyachi, and Tomita \cite{Ka_Mi_To2022} obtained a complete characterization in the case $\rho=0$ by working on local Hardy spaces $h^p$ and their endpoint counterparts $bmo$.

In the general $l$-linear case, for $0\le\delta\le \rho<1$ and exponents $0<p_1,\dots,p_l\le\infty$ with $\frac{1}{p_1}+\cdots+\frac{1}{p_l}=\frac{1}{p}$,
let
\[
m_0(\ppp)
   := -n\bigg(\sum_{j=1}^l \max\Bigl(\frac1{p_j},\frac12\Bigr)
        - \min\Bigl(\frac1p,\frac12\Bigr)\bigg),
\qquad
m_\rho(\ppp):=(1-\rho)\,m_0(\ppp),
\]
so that $m_0(\ppp)$ is the critical order in the case $\rho=0$, and $m_\rho(\ppp)$ is its natural $(1-\rho)$-dilation for $0<\rho<1$.  
As mentioned in \cite[Proposition 1.3]{Park_To2024}, combining the results of \cite{Ka_Mi_To2022} with a standard Littlewood-Paley decomposition and dilation argument, we can obtain the following unweighted estimate.
\begin{customthm}{E}\label{unweighted-pseudo-intro}\cite{Ka_Mi_To2022, Park_To2024}
Let $0\le\delta\le \rho<1$, $0<p_1,\dots,p_l\le\infty$, and  $\frac{1}{p_1}+\cdots+\frac{1}{p_l}=\frac{1}{p}$.  
Suppose that
\[
m < m_\rho(\ppp)
   =(1-\rho)\,m_0(\ppp),
\]
and 
\[
a\in \mathbb{M}_l S_{\rho,\delta}^m(\mathbb{R}^n).
\]
Then
\[
T_{[a]}:
H^{p_1}(\mathbb{R}^n)\times\cdots\times H^{p_l}(\mathbb{R}^n)\to X^p(\mathbb{R}^n)
\]
is bounded, where $H^p$ denotes the (real) Hardy space for $0<p<\infty$, we adopt the convention $H^\infty=L^\infty$, and
\[
X^p
  :=
  \begin{cases}
    L^p, & 0<p<\infty,\\[2pt]
    BMO, & p=\infty.
  \end{cases}
\]
In particular, if $1<p_1,\dots,p_l<\infty$, then
\[
T_{[a]}:
L^{p_1}(\mathbb{R}^n)\times\cdots\times L^{p_l}(\mathbb{R}^n)\to L^{p}(\mathbb{R}^n)
\]
is bounded under the same assumption $m<m_\rho(\ppp)$.
\end{customthm}

\medskip
Recently, the authors \cite{Park_To2024, Park_To2025} established weighted estimates for 
multilinear pseudo-differential operators in the multiple weight setting.

\begin{customthm}{F}\cite{Park_To2024, Park_To2025}\label{weightedmultire}
Let $0\le \delta\le \rho<1$, $1<r\le 2$, and $r<p_1,\dots,p_l <\infty$ with  $\frac{1}{p_1}+\cdots+\frac{1}{p_l}=\frac{1}{p}$. 
If
$a\in \mathbb{M}_l S_{\rho,\delta}^{-\frac{nl}{r}(1-\rho)}(\bbrn)$ and $\www=(w_1,\dots,w_l)\in A_{(\frac{p_1}{r},\dots,\frac{p_l}{r})}(\bbrn)$, then there exists a constant $C_{a,\www}>0$ such that
\[
\big\Vert T_{[a]}\big\Vert_{L^{p_1}(w_1)\times\cdots\times L^{p_l}(w_l)\to L^p(v_{\www})}
   \le C_{a,\www}.\]
\end{customthm} 
A fundamental question, analogous to the one raised in the multilinear Fourier multiplier setting, is whether the weight condition in Theorem \ref{weightedmultire} can be weakened.  
Our second main theorem shows that this is not possible.
\begin{theorem}\label{multipseudonega}
Let $l\in \bbn$, $0 \le \rho<1$, $1<r\le 2$, $0<q<\infty$, $\max\{r,q\} \le p_1,\dots,p_l<\infty$, and $\frac{1}{p_1}+\dots+\frac{1}{p_l}=\frac{1}{p}$.
Assume that for all  $a\in\mathbb{M}_lS_{\rho,0}^{-\frac{nl}{r}(1-\rho)}(\bbrn)$ and for all $\www=(w_1,\dots,w_l)\in A_{(\frac{p_1}{q},\dots,\frac{p_l}{q})}(\bbrn)$, there exists a constant $C_{a,\www}>0$ such that the weighted norm inequality
\begin{equation}\label{thm4maasu}
\big\Vert T_{[a]}\big\Vert_{L^{p_1}(w_1)\times\cdots\times L^{p_l}(w_l)\to L^p(v_{\www})}\le C_{a,\www} 
\end{equation}
holds.
Then it is necessary that
\[
r \le q.
\]
\end{theorem}
\begin{remark*}
In the proof of Theorem \ref{multipseudonega}, we will actually show that the condition $r \le q$ is necessary under a weaker assumption
on the symbol, namely that the weighted inequality \eqref{thm4maasu}
holds for all $x$-independent symbols $a$, and hence the above necessary
condition follows even in this restricted setting.
\end{remark*}

Using Rubio de Francia's extrapolation theorem, Chen \cite{Ch1991} obtained Theorem \ref{multipseudonega} with $l=1$. We will prove the multilinear case (including the linear case) without using the extrapolation theorem.

\medskip

In the same way as Chanillo and Torchinsky \cite{Ch_To1985},
Theorem \ref{weightedmultire} was proved by the following
sharp maximal function estimate for multilinear pseudo-differential operators (\cite{Park_To2024, Park_To2025}).
If $1<r \le 2$, $0 \le \delta \le \rho<1$, and $a\in \mathbb{M}_l S_{\rho,\delta}^{-\frac{nl}{r}(1-\rho)}(\bbrn)$,
then
\begin{equation}\label{ptmaxest}
\mathcal{M}^{\sharp}_{\frac{r}{l}}\big( T_{[a]}(f_1,\dots,f_l)\big)(x)\lesssim_a \mathbf{M}_r\big(f_1,\dots,f_l\big)(x), \q \q ~x\in\bbrn ,
\end{equation}
for all Schwartz functions $f_1, \dots f_l$ on $\bbrn$,
where
\[
\mathcal{M}^{\sharp}_{t}g(x)
=\sup_{Q \ni x}
\inf_{c_Q \in \mathbb{C}}
\bigg(\frac{1}{|Q|}\int_Q |g(y)-c_Q|^{t}\, dy\bigg)^{\frac{1}{t}}, \quad 0<t<\infty
\]
and $\mathbf{M}_r\big(f_1,\dots,f_l\big)(x)=\big( \mathbf{M}\big(|f_1|^r,\dots, |f_l|^r\big)(x)\big)^{\frac{1}{r}}$.
As a corollary of Theorem \ref{multipseudonega},
we obtain the optimality of the estimate \eqref{ptmaxest}.
\begin{corollary}\label{sharp-maximal-function-estimate}
Let $0 \le \rho<1$ and $1<r \le 2$.
If there exists $q>0$ such that
\begin{equation}\label{ass-smfe}
\mathcal{M}^{\sharp}_{\frac{r}{l}}\big( T_{[a]}(f_1,\dots,f_l)\big)(x)
\lesssim_a \mathbf{M}_q\big(f_1,\dots,f_l\big)(x), \q \q ~x\in\bbrn ,
\end{equation}
for all  $a\in \mathbb{M}_l S_{\rho,0}^{-\frac{nl}{r}(1-\rho)}(\bbrn)$,
then $r \le q$.
\end{corollary}
In fact, since Theorem \ref{multipseudonega} remains valid when restricted
to $x$-independent symbols, the above corollary holds even if
\eqref{ass-smfe} is assumed only for such symbols; see the proof of the corollary.

\bigskip

The proofs of Theorems~\ref{sharpnessmultilinear} and~\ref{multipseudonega} follow a common underlying idea.
We construct explicit families of multiple weights and model multipliers (or pseudo-differential symbols) that closely reflect the borderline behavior of the corresponding weighted estimates, and we compare upper and lower bounds obtained from suitably chosen test functions.

For Theorem \ref{sharpnessmultilinear}, we use the structural characterization of $A_{\ppp}$ in Lemma \ref{multiindiweight} to build a family of weights of the form
\[
\www=(w_{\beta_1,\gamma},\dots,w_{\beta_l,\gamma})
\qquad\text{with}\qquad
w_{\beta_j,\gamma}(x)=\frac{|x|^{\beta_j}}{|x-e_1|^{\gamma}}.
\]
These weights belong to $A_{(\frac{p_1}{q},\dots,\frac{p_l}{q})}(\bbrn)$ once the parameters $\beta_1,\dots,\beta_l$ and $\gamma$ are chosen appropriately. See Lemma \ref{mainweest} below.
We then combine them with prototype multipliers constructed from the Bessel potential kernel and apply the assumed inequality to localized bump functions.  
A comparison of exponents leads to the condition $q\ge \frac{nl}{s}$, which shows that the assumption is sharp.
The proof of Theorem \ref{multipseudonega} follows a similar strategy in the pseudo-differential setting.  
A Littlewood-Paley decomposition reduces the problem to symbols of type $(0,0)$, using interpolation together with the self-improving property of multiple weights.  
We then combine the same family of multiple weights with the pseudo-differential operator associated with $\sigma_{\frac{nl}{r}}$ and evaluate the assumed estimate on the previously introduced localized functions.  
The resulting comparison of exponents yields the condition $r\le q$, establishing the sharp necessity of the weighted pseudo-differential estimate.

\bigskip

\subsection*{Organization}
The paper is organized as follows.
In Section \ref{prelimsec} we collect basic properties of Muckenhoupt weights and multiple weights that will be used throughout the proofs.
 In Section \ref{keyestsec} we introduce the explicit families of multiple weights, multipliers, and test functions that underlie our sharpness arguments, and we establish the key weighted estimates for these examples. Section~4 is devoted to the proof of Theorem~\ref{sharpnessmultilinear}, where the constructions from Section~3 are used to derive the optimality of weight assumptions in the (multilinear) multiplier setting. Finally, in Section~5 we prove Theorem~\ref{multipseudonega} and Corollary \ref{sharp-maximal-function-estimate}.

\hfill

\section{Preliminaries}\label{prelimsec}

\subsection{Muckenhoupt's weight space $A_p(\mathbb{R}^n)$}

We recall several standard properties of the Muckenhoupt classes $A_p(\bbrn)$, which will be used later in the paper.  
First, for $1<p<\infty$, the duality relation
\begin{equation*}
w\in A_p(\bbrn)
\quad\Longleftrightarrow\quad
w^{1-p'}\in A_{p'}
\end{equation*}
is a direct consequence of the definition of the $A_p$ condition and H\"older's inequality.

Another useful fact is that a suitable combination of two $A_1$ weights yields an $A_p$ condition.  
If $w_0,w_1\in A_1(\bbrn)$, then
\begin{equation}\label{a1prodap}
w_0\, w_1^{\,1-p}\in A_p(\bbrn),
\end{equation}
that is, $A_1(\bbrn) A_1^{1-p}(\bbrn)\subset A_p(\bbrn)$.

As a basic model class, consider the power weights $w(x)= |x|^{\alpha}$.  
It is well known that such weights satisfy
\begin{equation}\label{singwexpower}
|x|^{\alpha} \in A_p(\bbrn)
\quad\Longleftrightarrow\quad
\begin{cases}
-n<\alpha<n(p-1),
&p>1,
\\
-n<\alpha \le 0,
&p=1,
\end{cases}
\end{equation}
which completely characterizes all radial power weights in $A_p(\bbrn)$.

We also recall from \cite[Chapter IV. Theorem 2.7]{Ga_Ru1985} the classical reverse H\"older self-improving property:
\begin{equation}\label{w1eapcon}
 w\in A_p(\mathbb{R}^n)
 \qquad\Longrightarrow\qquad
 w^{1+\epsilon}\in A_p(\mathbb{R}^n)
 \quad\text{for some }\epsilon>0.
\end{equation}

Finally, the scalar $A_p$ condition is invariant under dilations.  
For any $\lambda>0$,
\begin{equation}\label{ap-dilation}
w\in A_p(\mathbb{R}^n)
\quad\Longrightarrow\quad
w(\lambda \cdot )\in A_p(\mathbb{R}^n),
\end{equation}
and moreover, the $A_p$ constant of $w(\lambda\cdot)$ remains equal to that of $w$.  
This invariance property will be extended to multiple weights in the next subsection.

\subsection{Multiple weight space $A_{(p_1,\dots,p_l)}$}

We next summarize the structure and basic properties of the multiple weight classes $A_{\ppp}(\bbrn)$, which play a central role in the analysis of multilinear operators.  
In view of \eqref{multiplewedef}, the definition involves the combined weight
\[
v_{\www}:=\prod_{j=1}^l w_j^{\frac{p}{p_j}},
\qquad
\text{where}\quad 
\frac{1}{p}=\frac{1}{p_1}+\cdots+\frac{1}{p_l}.
\]
The following lemma provides a precise and complete characterization of $A_{\ppp}(\bbrn)$ in terms of classical scalar $A_p$ weights.
\begin{customlemma}{G}\cite[Theorem 3.6]{Le_Om_Pe_To_Tr2009}\label{multiindiweight}
Let $\www=(w_1,\dots,w_l)$ and $1\le p_1,\dots,p_l<\infty$. Then
\[
\www\in A_{\ppp}(\mathbb{R}^n)
\quad\Longleftrightarrow\quad
\begin{cases}
w_j^{\,1-p_j'}\in A_{lp_j'}(\mathbb{R}^n), & j=1,\dots,l,\\[2mm]
v_{\www}\in A_{lp}(\mathbb{R}^n),&
\end{cases}
\]
with the convention that when $p_j=1$, the condition $w_j^{\,1-p_j'}\in A_{lp_j'}(\bbrn)$ is replaced by $w_j^{\frac{1}{l}}\in A_1(\mathbb{R}^n)$.
\end{customlemma}

This structural decomposition shows that the multilinear class $A_{\ppp}(\bbrn)$ is governed entirely by scalar $A_p$ conditions.  
Thus, many classical properties of the $A_p(\bbrn)$ classes transfer directly to the multilinear setting once each component weight and the combined weight $v_{\www}$ are analyzed.  
In fact, Lemma~\ref{multiindiweight} is one of the key tools used in Chapter~3, where explicit multiple weights are constructed and their admissibility in $A_{\ppp}$ is verified by checking the corresponding scalar conditions.

The self-improving property \eqref{w1eapcon} of classical $A_p$ weights also holds in the multilinear setting.
Indeed, by Lemma~\ref{multiindiweight}, the condition $\www\in A_{\ppp}$ reduces to estimates for each $w_j$ and for the combined weight $v_{\www}$.
Applying \eqref{w1eapcon} to these weights yields a uniform $\epsilon>0$, we obtain the following lemma.
\begin{lemma}\cite[Lemma 4]{Park2024_submitted}\label{mwextra}
Let $1\le p_1,\dots,p_l<\infty$ and assume $\www\in A_{\ppp}(\mathbb{R}^n)$.
Then there exists $\epsilon>0$ such that
\[
\www^{1+\epsilon}
 :=
 (w_1^{1+\epsilon},\dots,w_l^{1+\epsilon})
 \in A_{\ppp}(\mathbb{R}^n).
\]
\end{lemma}

Another essential property of multiple weights is invariance under dilations, which is a multilinear analogue of \eqref{ap-dilation}.
Moreover, the invariance holds uniformly: the $A_{\ppp}$ constant does not depend on the dilation parameter.

\begin{lemma}\label{dilaaplemma}
For any $\www\in A_{\ppp}(\mathbb{R}^n)$,
\[
  \www(\lambda\cdot)
  :=
  (w_1(\lambda\cdot),\dots,w_l(\lambda\cdot))
  \in A_{\ppp}(\mathbb{R}^n)
\]
uniformly in $\lambda>0$; that is, the $A_{\ppp}$ constant of $\www(\lambda\cdot)$ can be chosen
independently of $\lambda$.
\end{lemma}
The lemma follows directly from
\[
[\www(\lambda\cdot)]_{A_{\ppp}}= [\www]_{A_{\ppp}},
\]
which is an immediate consequence of the definition \eqref{multiplewedef}.

\subsection{Bessel potential}\label{besselpre}

Fix \(N\in\mathbb{N}\) and consider functions on \(\mathbb{R}^N\).  
For \(t>0\), the Bessel potential kernel \(G_t\) is defined as the Fourier transform of
\[
(1+4\pi^2 |x|^2)^{-\frac{t}{2}}, \qquad x\in\mathbb{R}^N.
\]

In this paper we only require the behavior of \(G_t\) near the origin in the range \(0<t \le N\), 
for which
\begin{equation}\label{bessel2}
G_t(\xi)\sim_{t,N}\begin{cases} |\xi|^{-(N-t)}, &0<t<N\\
\ln\big(2|\xi|^{-1}\big),&t=N
\end{cases}
\end{equation}
when $|\xi|\le 1$.
Thus, for \(0<t<N\), the kernel \(G_t\) behaves like the homogeneous function
\(|\xi|^{-(N-t)}\) near the origin, while in the critical case \(t=N\) it exhibits
a logarithmic singularity.
When \(t>N\), the kernel \(G_t\)  has no singularity at the origin.
We refer to \cite[Chapter~1.2.2]{Gr2} for more details.

\hfill

\section{Key Estimates}\label{keyestsec}
We fix the multilinearity $l\in\bbn$ and construct examples of weights, operators, and functions, which will be employed in the proof of Theorems \ref{sharpnessmultilinear} and \ref{multipseudonega}.

\subsection{Contruction of (multiple) weights}
Let $e_1=(1,0,\dots,0)\in \bbrn$.
Let $0<p_1,\dots,p_l<\infty$,
$\frac{1}{p_1}+\dots+\frac{1}{p_l}=\frac{1}{p}$, and $0<q\le \min\big\{p_1,\dots,p_l\big\}$.
For any $\beta_1,\dots,\beta_l\ge 0$ and $\gamma>0$, we define
\begin{equation}\label{wbegade}
w_{\beta_j,\gamma}(x)=\frac{|x|^{\beta_j}}{|x-e_1|^{\gamma}},\quad x\in\bbrn.
\end{equation}
\begin{lemma}\label{mainweest}
Suppose that $0<\gamma<n$, and for each $j$ we set $0<\beta_j<n(\frac{p_j}{q}-1)$ when $p_j>q$, while we assume $\beta_j=0$ in the case $p_j=q$.
Then we have
\begin{equation}\label{multsharpkeyest}
\www:=(w_{\beta_1,\gamma},\dots,w_{\beta_l,\gamma})\in A_{(\frac{p_1}{q},\dots,\frac{p_l}{q})}(\bbrn).
\end{equation}
\end{lemma}
\begin{proof}

According to Lemma \ref{multiindiweight},
the assertion \eqref{multsharpkeyest} is equivalent to
\begin{equation}\label{sepawecon}
\begin{cases}
\big( w_{\beta_j,\gamma}\big)^{1-(\frac{p_j}{q})'}\in A_{l(\frac{p_j}{q})'}(\bbrn), & j=1,\dots,l, \\
v_{\www}\in A_{\frac{lp}{q}}(\bbrn),
\end{cases}
\end{equation}
where the condition $\big(w_{\beta_j,\gamma}\big)^{1-(\frac{p_j}{q})'}\in A_{l(\frac{p_j}{q})'}(\bbrn)$ in the case $p_j=q$ is understood as $\big(w_{\beta_j,\gamma}\big)^{\frac{1}{l}}=\big(w_{0,\gamma}\big)^{\frac{1}{l}}\in A_1(\bbrn)$.

If $q<p_j$, then
$$0<\beta_j\Big(\big(\frac{p_j}{q}\big)'-1\Big)<n \quad \text{ and } \quad 0<\gamma\frac{\big(\frac{p_j}{q}\big)'-1}{l\big(\frac{p_j}{q}\big)'-1}\le \gamma<n,$$
and thus \eqref{singwexpower} yields that
$$ |x|^{-\beta_j((\frac{p_j}{q})'-1)}\in A_1(\bbrn)      \quad \text{ and } \quad |x-e_1|^{-\gamma \frac{(\frac{p_j}{q})'-1}{l(\frac{p_j}{q})'-1}}\in A_1(\bbrn).$$
This deduces, in view of \eqref{a1prodap},
$$\big(w_{\beta_j,\gamma}(x)\big)^{1-(\frac{p_j}{q})'}=|x|^{-\beta_j((\frac{p_j}{q})'-1)} \Big(|x-e_1|^{-\gamma \frac{(\frac{p_j}{q})'-1}{l(\frac{p_j}{q})'-1}} \Big)^{1-l(\frac{p_j}{q})'}\in A_{l(\frac{p_j}{q})'}(\bbrn).$$
If $q=p_j$, then
$$\big(w_{\beta_j,\gamma}(x)\big)^{\frac{1}{l}}=\big(w_{0,\gamma}(x)\big)^{\frac{1}{l}}=|x-e_1|^{-\frac{\gamma}{l}}\in A_1(\bbrn)$$
because $0<\frac{\gamma}{l}\le \gamma<n$.
This concludes the first condition of \eqref{sepawecon}.

To verify the second one of \eqref{sepawecon},
we write
\begin{equation}\label{wwxdecom}
v_{\www}(x)=\frac{|x|^{\frac{p}{p_1}\beta_1+\dots+\frac{p}{p_l}\beta_l}}{|x-e_1|^{\gamma}}.
\end{equation}
Clearly, in view of \eqref{singwexpower},
\begin{equation}\label{gateina1}
\frac{1}{|x-e_1|^{\gamma}}\in A_1(\bbrn).
\end{equation}
If $p_1=\dots=p_l=q$, then
$$v_{\www}(x)=\frac{1}{|x-e_1|^{\gamma}}\in A_1(\bbrn)=A_{\frac{lp}{q}}(\bbrn).$$
If $q<\max\{p_1,\dots,p_l\}$, then
$$0<\frac{p}{p_1}\beta_1+\dots+\frac{p}{p_l}\beta_l<  \sum_{j=1}^{l}\frac{pn}{p_j}(\frac{p_j}{q}-1)=n\Big(\frac{lp}{q}-1\Big),$$
which implies
$$0<\frac{\frac{p}{p_1}\beta_1+\dots+\frac{p}{p_l}\beta_l}{\frac{lp}{q}-1}<n$$
so that
\begin{equation}\label{xmppjbjina1}
|x|^{-\frac{\frac{p}{p_1}\beta_1+\dots+\frac{p}{p_l}\beta_l}{\frac{lp}{q}-1}}\in A_1(\bbrn).
\end{equation}
Combining \eqref{wwxdecom}, \eqref{gateina1}, and \eqref{xmppjbjina1}, the property \eqref{a1prodap} deduces that
$$v_{\www}(x)=\frac{1}{|x-e_1|^{\gamma}}\Big( |x|^{-\frac{\frac{p}{p_1}\beta_1+\dots+\frac{p}{p_l}\beta_l}{\frac{lp}{q}-1}}\Big)^{1-\frac{lp}{q}}\in A_{\frac{lp}{q}}(\bbrn),$$
which completes the proof of \eqref{sepawecon} and thus \eqref{multsharpkeyest} also holds.
\end{proof}

\subsection{Construction of (multilinear) multipliers}
For $0<\mu\le nl$ we define a function $\sigma_{\mu}$ on $(\bbrn)^l$ as
\begin{equation*}
\sigma_{\mu}(\xxxi)=\frac{1}{(1+4\pi^2|\xxxi|^2)^{\frac{\mu}{2}}}e^{-2\pi i\langle e_1,\xi_1+\dots+\xi_l\rangle}
\end{equation*}
where $\xxxi=(\xi_1,\dots,\xi_l)\in (\bbrn)^l$.
We first claim that 
\begin{equation}\label{mulkeyest}
\sup_{k\in\bbz}\Big\Vert \sigma_{\mu}(2^k\vec{\cdot}\,)\wh{\Psi}\Big\Vert_{L^2_{\mu}(\bbr^{nl})}\lesssim 1.
\end{equation}
Indeed, choosing a positive integer $M$ with $0<\mu<2M$, and using Plancherel's theorem and H\"older's inequality with $\frac{2M}{\mu}>1$, we obtain
\begin{align*}
\Big\Vert \sigma_{\mu}(2^k\vec{\cdot}\,)\wh{\Psi}\Big\Vert_{L^2_{\mu}(\bbr^{nl})}&=\Big\Vert   \big(1+4\pi^2|\vec{\cdot}\,|^2 \big)^{\frac{\mu}{2}}\Big(\sigma_{\mu}(2^k\vec{\cdot}\,)\wh{\Psi} \Big)^{\vee}   \Big\Vert_{L^2(\bbr^{nl})}\\
&\le \Big\Vert  \Big(\sigma_{\mu}(2^k\vec{\cdot}\,)\wh{\Psi} \Big)^{\vee}\Big\Vert_{L^2(\bbr^{nl})}^{1-\frac{\mu}{2M}}\Big\Vert \big( 1+4\pi^2|\vec{\cdot}\,|^2\big)^{M}\Big( \sigma_{\mu}(2^k\vec{\cdot}\,)\wh{\Psi}\Big)^{\vee}\Big\Vert_{L^2(\bbr^{nl})}^{\frac{\mu}{2M}}\\
&=\Big\Vert \sigma_{\mu}(2^k\vec{\cdot}\,)\wh{\Psi}\Big\Vert_{L^2(\bbr^{nl})}^{1-\frac{\mu}{2M}}\Big\Vert \big( \vec{I}-\vec{\Delta}\big)^{M}\big( \sigma_{\mu}(2^k\vec{\cdot}\,)\wh{\Psi}\big)\Big\Vert_{L^2(\bbr^{nl})}^{\frac{\mu}{2M}}.
\end{align*}
Since
\begin{align*}
\Big| \sigma_{\mu}(2^k\xxxi)\wh{\Psi}(\xxxi)\Big|\sim 2^{-\mu k}\chi_{|\xxxi|\sim 1},
\end{align*}
we have
$$\Big\Vert \sigma_{\mu}(2^k\vec{\cdot}\,)\wh{\Psi}\Big\Vert_{L^2(\bbr^{nl})}\sim 2^{-\mu k}.$$
Similarly, a straightforward computation leads to
\begin{align*}
\Big| \big( \vec{I}-\vec{\Delta}\big)^M\big(\sigma_{\mu}(2^k\vec{\cdot}\, )\wh{\Psi}\big)(\xxxi)\Big|\sim 2^{-\mu k}2^{2Mk}\chi_{|\xxxi|\sim 1},
\end{align*}
which implies that
$$\Big\Vert \big( \vec{I}-\vec{\Delta}\big)^{M}\big( \sigma_{\mu}(2^k\vec{\cdot}\,)\wh{\Psi}\big)\Big\Vert_{L^2(\bbr^{nl})}\sim 2^{-\mu k}2^{2Mk}.$$
This yields that
\begin{equation*}
\Big\Vert \sigma_{\mu}(2^k\vec{\cdot}\,)\wh{\Psi}\Big\Vert_{L^2_{\mu}(\bbr^{nl})}\lesssim 2^{-\mu k(1-\frac{\mu}{2M})}\big( 2^{-\mu k}2^{2Mk}\big)^{\frac{\mu}{2M}}=1,
\end{equation*}
which completes the proof of \eqref{mulkeyest}.

In addition, a straightforward calculation shows that 
\begin{equation}\label{multipseudoch}
\sigma_{\mu}\in \mathbb{M}^lS_{(0,0)}^{-\mu}(\bbrn).
\end{equation}
Then we consider the corresponding (multilinear) multiplier operator $T_{\sigma_{\mu}}$ as
 $$T_{\sigma_{\mu}}\big(f_1,\dots,f_l\big)(x)=\int_{(\bbrn)^l}\sigma_{\mu}(\xi_1,\dots,\xi_l)\Big( \prod_{j=1}^{l}\wh{f_j}(\xi_j)\Big) \,e^{2\pi i\langle x,\xi_1+\dots+\xi_l\rangle} \; d\xxxi$$
for $f_1,\dots,f_l\in \mathscr{S}(\bbrn)$. 
In this case, the operator $T_{\sigma_{\mu}}$ can be written as 
\begin{align*}
T_{\sigma_{\mu}}\big(f_1,\dots,f_l \big)(x)&=\int_{(\bbrn)^l}G_\mu\big(x-e_1-y_1,\dots,x-e_1-y_l \big)\prod_{j=1}^{l}f_j(y_j)\, d\yyy
\end{align*}
where we recall $G_{\mu}$ is the Bessel potential kernel on $\bbr^{nl}$, discussed in Section \ref{besselpre}.

\subsection{Main weighted estimates}
Let $\phi$ be a Schwartz function on $\bbrn$ such that
$\supp(\phi)\subset \{x\in\bbrn: |x| \le 2 \}$, $0\le \phi\le 1$, and $\phi(x)=1$ for all $x$ with $|x|\le 1$.
For arbitrary $0<\epsilon<\frac{1}{4l}$, we define a function $f^{\epsilon}$ on $\bbrn$ to be
$$f^{\epsilon}(x)=\phi(\epsilon^{-1}x)$$
so that $0\le f^{\epsilon}\le 1$ and
$$f^{\epsilon}(x)=\begin{cases}  1& |x|\le \epsilon \\ 0 &|x|\ge 2\epsilon\end{cases}.$$
Setting
\begin{equation}\label{fjdef}
f_1(x)=\dots=f_l(x)=f^{\epsilon}(x)=\phi(\epsilon^{-1}x),
\end{equation}
we have
\begin{align*}
T_{\sigma_{\mu}}\big(f^{\epsilon},\dots,f^{\epsilon} \big)(x)&=\int_{(\bbrn)^l}G_{\mu}\big(x-e_1-y_1,\dots,x-e_1-y_l \big)\prod_{j=1}^{l}f^{\epsilon}(y_j)\, d\yyy\\
&\ge \int_{|y_1|,\dots,|y_l|<\epsilon}G_{\mu}\big(x-e_1-y_1,\dots,x-e_1-y_l \big)\, d\yyy.
\end{align*}
If $|y_1|,\dots, |y_l|<\epsilon$ and $|x-e_1|<3\epsilon$, then
$$\big| (x-e_1-y_1,\dots,x-e_1-y_l)\big|\le  \sum_{j=1}^{l}|x-e_1-y_j|<4l\epsilon<1,$$
for which
$$G_{\mu}(x-e_1-y_1,\dots,x-e_1-y_l)\gtrsim \begin{cases}\frac{1}{\epsilon^{nl-\mu}} & \mu<nl\\ \ln(\frac{1}{2l\epsilon})& \mu=nl \end{cases},$$
in view of \eqref{bessel2}.
This yields that for $|y_1|,\dots, |y_l|<\epsilon$ and $|x-e_1|<3\epsilon$,
$$G_{\mu}(x-e_1-y_1,\dots,x-e_1-y_l)\gtrsim \frac{1}{\epsilon^{nl-\mu}}$$
and thus
\begin{equation}\label{multiptlowerest}
T_{\sigma_{\mu}}\big(f^{\epsilon},\dots,f^{\epsilon} \big)(x)\gtrsim \epsilon^{\mu} \q \text{ for }~ |x-e_1|<3\epsilon.
\end{equation}
Let $\www$ be the (multiple) weights in \eqref{multsharpkeyest}.
Then for $2\epsilon<|x-e_1|<3\epsilon$,
$$v_{\www}(x)\sim \frac{1}{|x-e_1|^{\gamma}}\sim\epsilon^{-\gamma}$$
and thus, employing \eqref{multiptlowerest}, we obtain
\begin{equation}\label{lowerlpvwest}
\big\Vert T_{\sigma_{\mu}}(f^{\epsilon},\dots,f^{\epsilon})\big\Vert_{L^p(v_{\www})}\ge \bigg(  \int_{2\epsilon<|x-e_1|<3\epsilon} \big| T_{\sigma_{\mu}}\big(f^{\epsilon},\dots, f^{\epsilon}\big)(x)\big|^p v_{\www}(x)
\, dx \bigg)^{\frac{1}{p}}\gtrsim \epsilon^{\mu}\epsilon^{-\frac{\gamma}{p}}\epsilon^{\frac{n}{p}}.
\end{equation}\label{operlowerest}
In addition, for each $j=1,\dots,m$,
if $|x|\le 2\epsilon$, then 
$$w_{\beta_j,\gamma}(x)\sim_{\gamma} |x|^{\beta_j}\lesssim \epsilon^{\beta_j},$$
and thus
\begin{equation}\label{feplpjest}
\Vert f^{\epsilon} \Vert_{L^{p_j}(w_{\beta_j,\gamma})}\le \bigg( \int_{|x|\le 2\epsilon}w_{\beta_j,\gamma}(x)\, dx\bigg)^{\frac{1}{p_j}}\lesssim \epsilon^{\frac{\beta_j}{p_j}}\epsilon^{\frac{n}{p_j}}.
\end{equation}

\hfill

\section{Proof of Theorem \ref{sharpnessmultilinear}}

We assume that $q<\frac{nl}{s}$ and will derive a contradiction.
Choose a constant $\delta>0$ such that
\begin{equation}\label{condemult}
s+\frac{2\delta}{p}<\frac{nl}{q} \q \text{ and }\q \delta<n\min_{1 \le j \le l}\Big\{\Big(\frac{p_j}{q}-1 \Big),1 \Big\}.
\end{equation}
Take $$\gamma=n-\delta>0$$
and $\beta_1,\dots,\beta_l\ge 0$ to be
\begin{equation}\label{betajdef}
\beta_j= \begin{cases} 0 & ~\text{ if }~ p_j=q,
\\
n(\frac{p_j}{q}-1)-\delta
& ~\text{ if }~ p_j>q, \end{cases}
\end{equation}
and then we set
$$\www=(w_{\beta_1,\gamma},\dots,w_{\beta_l,\gamma}).$$ 
Since the parameters $\gamma,\beta_1,\dots,\beta_l$ satisfy the assumptions in Lemma \ref{mainweest},
we have
\begin{equation}\label{wwwinapjq}
\www\in A_{(\frac{p_1}{q},\dots,\frac{p_l}{q})}(\bbrn).
\end{equation}

Now we set
\begin{equation*}
\sigma(\xxxi)=\sigma_{s}(\xxxi), \quad \xxxi=(\xi_1,\dots,\xi_l)\in (\bbrn)^l.
\end{equation*}
According to \eqref{mulkeyest},
\begin{equation}\label{supksi2kest}
\sup_{k\in\bbz}\big\Vert \sigma(2^k\vec{\cdot}\,)\wh{\Psi}\big\Vert_{L^2_s(\bbr^{nl})}\lesssim 1.
\end{equation}
For any $0<\epsilon<\frac{1}{4l}$, we set
\begin{equation*}
f_1(x)=\dots=f_l(x)=f^{\epsilon}(x)
\end{equation*}
as in \eqref{fjdef}.
Then it follows from \eqref{lowerlpvwest}, \eqref{multiweibound}, \eqref{supksi2kest}, and  \eqref{feplpjest} that
$$\epsilon^{s}\epsilon^{-\frac{\gamma}{p}}\epsilon^{\frac{n}{p}}\lesssim \big\Vert T_{\sigma}(f^{\epsilon},\dots,f^{\epsilon})\big\Vert_{L^p(v_{\www})}\lesssim \prod_{j=1}^{l}\Vert f^{\epsilon}\Vert_{L^{p_j}(w_{\beta_j,\gamma})}\lesssim \epsilon^{\frac{\beta_1}{p_1}+\dots+\frac{\beta_l}{p_l}}\epsilon^{\frac{n}{p}}.$$
Thus, since $\frac{\beta_1}{p_1}+\dots+\frac{\beta_l}{p_l} \ge \frac{nl}{q}-\frac{n}{p}-\frac{\delta}{p}$,
$$1\lesssim \epsilon^{\frac{\beta_1}{p_1}+\dots+\frac{\beta_l}{p_l}}\epsilon^{-s}\epsilon^{\frac{\gamma}{p}}\le \epsilon^{\frac{nl}{q}-s-\frac{2\delta}{p}},$$
which clearly leads to a contradiction as $\epsilon$ goes to $0$, due to the first one in \eqref{condemult}.

This completes the proof.

\hfill

\section{ Proof of Theorem \ref{multipseudonega} and Corollary \ref{sharp-maximal-function-estimate}}

\subsection{Reduction to symbols of type $(0,0)$}
Assume that \eqref{thm4maasu} holds for all $x$-independent symbols $a\in\mathbb{M}_lS_{\rho,0}^{-\frac{nl(1-\rho)}{r}}(\bbrn)$ and all (multiple) weights $\www\in A_{(\frac{p_1}{q},\dots,\frac{p_l}{q})}(\bbrn)$. Then we  claim that for all $x$-independent symbols $b\in\mathbb{M}_lS_{0,0}^{-\frac{nl}{r}}(\bbrn)$ and all (multiple) weights $\www\in A_{(\frac{p_1}{q},\dots,\frac{p_l}{q})}(\bbrn)$, there exists a constant $C_{b,\www}>0$ such that
\begin{equation}\label{claimtbpseu}
\big\Vert T_{[b]}\big\Vert_{L^{p_1}(w_1)\times\cdots\times L^{p_l}(w_l)\to L^p(v_{\www})}\le C_{b, \www}
\end{equation}
holds.
Since it is exactly the estimate \eqref{thm4maasu} when $\rho=0$, we only need to consider the case $0<\rho<1$ for the claim.
Let $\{\Phi_j\}_{j\in\bbn_0}$ be a family of (inhomogeneous) Littlewood-Paley functions on $(\bbrn)^l$ defined with respect to the $\rho$-dependent dyadic scaling, so that
 $\supp(\wh{\Phi_0})\subset \{\xxxi\in (\bbrn)^l:|\xxxi|\le 2\}$, $\supp(\wh{\Phi_k})\subset \{\xxxi\in (\bbrn)^l: 2^{k-1}\le |\xxxi|\le 2^{k+1}\}$, and $\sum_{k=0}^{\infty}\wh{\Phi_k}(2^{k\rho}\xxxi)=1$.
Such a system can be obtained in the usual manner: choose a smooth cutoff $\wh{\Phi_0}$ satisfying $\wh{\Phi_0}(\xxxi)=1$ for $|\xxxi|\le 1$ and $\supp(\wh{\Phi_0})\subset \{\xxxi\in (\bbrn)^l:|\xxxi|\le 2 \}$, and define $\widehat{\Phi}(\xxxi)=\widehat{\Phi_0}(\xxxi)-\widehat{\Phi_0}(2^{1-\rho}\xxxi)$ and $\wh{\Phi_k}(\xxxi)=\wh{\Phi}(2^{-k}\xxxi)$ for $k\ge 1$.
For each $k\in\bbn_0$, let $$\lambda_k(\xxxi):=\wh{\Phi_k}(2^{k\rho}\xxxi),$$
and then we define
$b_k(\xxxi):=\lambda_k(\xxxi)b(\xxxi)$ so that
$$b(\xxxi)=\sum_{k=0}^{\infty}b_k(\xxxi).$$
Now let 
$$a_k(\xxxi):=b_k(2^{-k\rho}\xxxi)=\lambda_k(2^{-k\rho}\xxxi)b(2^{-k\rho}\xxxi)= \wh{\Phi_k}(\xxxi)b(2^{-k\rho}\xxxi).$$
Then we have
$$\supp(a_k)\subset \big\{ \xxxi\in (\bbrn)^l: 2^{k-1}\le |\xxxi|\le 2^{k+1}\big\},$$
where $\{2^{k-1} \le |\xxxi| \le 2^{k+1}\}$ is replaced by $\{|\xxxi| \le 2\}$ if $k=0$,
and  for any multi-indices $\vec{\alpha}\in (\bbn_0^n)^l$, there exists a constant $C_{b,\vec{\alpha}}$, independent of $k$, such that
$$\big| \partial_{\xxxi}^{\vec{\alpha}}a_k(\xxxi)\big|\le C_{b,\vec{\alpha}} \big(1+|\xxxi|\big)^{-\frac{nl(1-\rho)}{r}-\rho|\vec{\alpha}|},$$
that is,
$$a_k\in\mathbb{M}_lS_{\rho,0}^{-\frac{nl(1-\rho)}{r}}(\bbrn) \quad \text{ uniformly in }~k\in\bbn_0.$$

Assume  $\www\in A_{(\frac{p_1}{q},\dots,\frac{p_l}{q})}(\bbrn)$ and applying Lemma \ref{mwextra} to take $\epsilon>0$ so that $\www^{1+\epsilon}=(w_1^{1+\epsilon},\dots,w_l^{1+\epsilon})\in A_{(\frac{p_1}{q},\dots,\frac{p_l}{q})}(\bbrn)$. 
By the boundedness assumption  \eqref{thm4maasu},
$$\big\Vert T_{[a_k]}(f_1,\dots,f_l)\big\Vert_{L^p(v_{\www^{1+\epsilon}})}
\lesssim_{b, \www} \prod_{j=1}^{l}\Vert f_j\Vert_{L^{p_j}(w_j^{1+\epsilon})}.$$
Since $r \le \min\{p_1,\dots,p_l,2\}$, we see that
\[
\sum_{j=1}^l \max\left(\frac{1}{p_j},\frac{1}{2}\right)-\min\left(\frac{1}{p},\frac{1}{2}\right)
<\frac{l}{r},
\]
which implies $-\frac{nl(1-\rho)}{r}<m_{\rho}(\ppp)$.
Consequently,  Theorem \ref{unweighted-pseudo-intro} yields that
$$\big\Vert T_{[a_k]}(f_1,\dots,f_l)\big\Vert_{L^p(\bbrn)}\lesssim_b 2^{-\epsilon_0 k}\prod_{j=1}^{l}\Vert f_j\Vert_{L^{p_j}(\bbrn)}$$
for some $\epsilon_0>0$, depending on $\rho,p_1,\dots,p_l$.
Then we apply complex interpolation in \cite[Theorem 3.1]{COY2022} to obtain
\begin{equation}\label{ajreduction}
\big\Vert T_{[a_k]}(f_1,\dots,f_l)\big\Vert_{L^p(v_{\www})}
\lesssim_{b,\www} 2^{-\epsilon_0 \frac{\epsilon}{1+\epsilon}k}\prod_{j=1}^{l}\Vert f_j\Vert_{L^{p_j}(w_j)}.
\end{equation}

Now we write
$$T_{[b]}\big(f_1,\dots,f_l\big)=\sum_{k=0}^{\infty}T_{[b_k]}\big(f_1,\dots,f_l\big).$$
For each $k\in\bbn_0$, setting $f_j^{\rho}(x):=f_j(2^{k\rho}x)$ for $1\le j\le l$, we have
\begin{align*}
T_{[b_k]}\big(f_1,\dots,f_l\big)(x)&=\int_{(\bbrn)^l}b_k(\xxxi)\Big( \prod_{j=1}^{l}\wh{f_j}(\xi_j)\Big) \, e^{2\pi i\langle x,\xi_1+\dots+\xi_l\rangle} \; d\xxxi\\
&=\int_{(\bbrn)^l}a_k(\xxxi)\bigg( \prod_{j=1}^{l}\Big( 2^{-k\rho n}\wh{f_j}(2^{-k\rho}\xi_j)\Big)\bigg) e^{2\pi i\langle 2^{-k\rho}x,\xi_1+\dots+\xi_l\rangle} \; d\xxxi\\
&=T_{[a_k]}\big(f_1^{\rho},\dots,f_l^{\rho} \big)(2^{-k\rho}x).
\end{align*}

Note that 
$\big( w_1(2^{k\rho}{\cdot}),\dots,    w_l(2^{k\rho}{\cdot})      \big)\in A_{(\frac{p_1}{q},\dots,\frac{p_l}{q})}(\bbrn)$ uniformly in $k\in\mathbb{N}_0$ (see Lemma \ref{dilaaplemma})
and then the estimate \eqref{ajreduction} yields that
\begin{align*}
\big\Vert T_{[b_k]}\big(f_1,\dots,f_l\big)\big\Vert_{L^p(v_{\www})}&=2^{\frac{k\rho n}{p}}\Big\Vert T_{[a_k]}\big(f_1^{\rho},\dots,f_l^{\rho}\big)\Big\Vert_{L^p(v_{\www}(2^{k\rho}\cdot))}\\
&\lesssim_{b,\www} 2^{\frac{k\rho n}{p}} 2^{-\epsilon_0 \frac{\epsilon}{1+\epsilon}k}\prod_{j=1}^{l}\big\Vert f_j^{\rho}\big\Vert_{L^{p_j}(w_j(2^{k\rho}\cdot))}\\
&= 2^{-\epsilon_0 \frac{\epsilon}{1+\epsilon}k}\prod_{j=1}^{l}\big\Vert f_j\big\Vert_{L^{p_j}(w_j)}.
\end{align*}
This deduces that
\begin{align*}
\big\Vert T_{[b]}(f_1,\dots,f_l)\big\Vert_{L^p(v_{\www})}&\le \bigg(\sum_{k=0}^{\infty}\big\Vert T_{[b_k]}\big(f_1,\dots,f_l\big)\big\Vert_{L^p(v_{\www})}^{\min\{1,p\}} \bigg)^{\frac{1}{\min\{1,p\}}}\lesssim_{b,\www} \prod_{j=1}^{l}\big\Vert f_j\big\Vert_{L^{p_j}(w_j)},
\end{align*}
which completes the proof of \eqref{claimtbpseu}.

\subsection{Proof of Theorem \ref{multipseudonega} for symbols of type $(0,0)$.}
Due to the reduction carried out above, we may assume without loss of generality that $\rho=0$.
Specifically, suppose that  for all $x$-independent symbols $b\in\mathbb{M}_lS_{0,0}^{-\frac{nl}{r}}(\bbrn)$
the estimate \eqref{claimtbpseu} holds where $\www\in A_{(\frac{p_1}{q},\dots,\frac{p_l}{q})}(\bbrn)$. Under this assumption, we will show that necessarily $r\le q$.
To this end, assume towards contradiction that $r>q$.
We first choose $\delta>0$ such that
$$\frac{nl}{r}+\frac{2\delta}{p}<\frac{nl}{q}\quad \text{and}\quad \delta<n\min_{1\le j\le l} \Big\{\Big(\frac{p_j}{q}-1 \Big), 1 \Big\}.$$
and take $$\gamma=n-\delta>0$$
and $\beta_1,\dots,\beta_l>0$ satisfying \eqref{betajdef}.
Then it was already verified in \eqref{wwwinapjq} that 
$$\www=(w_{\beta_1, \gamma},\dots,w_{\beta_l,\gamma})\in A_{(\frac{p_1}{q},\dots,\frac{p_l}{q})}(\bbrn)$$
where $w_{\beta_j,\gamma}$ is defined as in \eqref{wbegade}.

Letting $0<\epsilon<\frac{1}{4l}$ and setting
$b(x,\xxxi)=\sigma_{\frac{nl}{r}}(\xxxi)$
and
\begin{equation*}
f_1(x)=\dots=f_l(x)=f^{\epsilon}(x)= \phi(\epsilon^{-1}x),
\end{equation*}
we may employ \eqref{multipseudoch}, \eqref{lowerlpvwest}, and \eqref{feplpjest} as 
$$T_{[b]}\big(f_1,\dots,f_l\big)=T_{\sigma_{\frac{nl}{r}}}\big(f^{\epsilon},\dots,f^{\epsilon}\big).$$
In other words,
$$b\in \mathbb{M}S_{(0,0)}^{-\frac{nl}{r}}(\bbrn)$$
and the estimate \eqref{claimtbpseu} yields
\begin{align*}
 \epsilon^{\frac{nl}{r}}\epsilon^{-\frac{\gamma}{p}}\epsilon^{\frac{n}{p}}\lesssim \Big\Vert T_{[b]}\big(f^{\epsilon},\dots,f^{\epsilon} \big) \Big\Vert_{L^p(v_{\ppp})}\lesssim \prod_{j=1}^{l}\Vert f^{\epsilon}\Vert_{L^{p_j}(w_j)}\lesssim     \epsilon^{\frac{\beta_1}{p_1}+\dots+\frac{\beta_l}{p_l}}\epsilon^{\frac{n}{p}},
\end{align*}
which deduces
$$1\lesssim  \epsilon^{\frac{\beta_1}{p_1}+\dots+\frac{\beta_l}{p_l}+\frac{\gamma}{p}-\frac{nl}{r}}\le \epsilon^{\frac{nl}{q}-\frac{nl}{r}-\frac{2\delta}{p}}$$
for all $0<\epsilon<\frac{1}{4l}$.
This is a contradiction when $\epsilon>0$ is arbitrary small.

This finishes the proof.

\subsection{Proof of Corollary \ref{sharp-maximal-function-estimate}}

Let $\rho, r$ be as in the assumption of Corollary \ref{sharp-maximal-function-estimate},
and assume that \eqref{ass-smfe} holds for some $q>0$.
We take $r,q <p_1,\dots,p_l<\infty$ and $0<p<\infty$ satisfying $\frac{1}{p_1}+\dots+\frac{1}{p_l}=\frac{1}{p}$.
If we can prove that \eqref{thm4maasu} holds for
all $x$-independent symbols $a\in\mathbb{M}_lS_{\rho,0}^{-\frac{nl}{r}(1-\rho)}(\bbrn)$ and
$\www=(w_1,\dots,w_l)\in A_{(\frac{p_1}{q},\dots,\frac{p_l}{q})}(\bbrn)$,
then $r \le q$ by Theorem  \ref{multipseudonega}.

Let $a$ be an $x$-independent symbol in
$\mathbb{M}_lS_{\rho,0}^{-\frac{nl}{r}(1-\rho)}(\bbrn)$,
$\www=(w_1,\dots,w_l)\in A_{(\frac{p_1}{q},\dots,\frac{p_l}{q})}(\bbrn)$,
and $f_1,\dots,f_l \in \mathscr{S}(\bbrn)$.
It follows from Lemma \ref{multiindiweight} that $v_{\www}\in A_{\frac{lp}{q}}(\mathbb{R}^n)$.
We take $0<t < \min\left\{1, p, \frac{r}{l}, \frac{q}{l}\right\}$.
Note that $1<\frac{p}{t}<\infty$ and $\frac{lp}{q}<\frac{p}{t}$,
and the latter gives $v_{\www} \in A_{\frac{p}{t}}$ by \eqref{apweightinclusion}.
Hence, the Hardy-Littlewood maximal operator $\mathcal{M}$ is bounded on $L^{\frac{p}{t}}(v_{\www})$,
and consequently $\mathcal{M}\big(|T_{[a]}(f_1,\dots,f_l)|^t\big) \in L^{\frac{p}{t}}(v_{\www})$.
Then
\begin{align*}
&\big\|T_{[a]}(f_1,\dots,f_l)\big\|_{L^p(v_{\www})}
=\left\| \big|T_{[a]}(f_1,\dots,f_l)\big|^t \right\|_{L^{\frac{p}{t}}(v_{\www})}^{\frac{1}{t}}
\le \left\|\mathcal{M}\left(\big|T_{[a]}(f_1,\dots,f_l)\big|^t\right) \right\|_{L^{\frac{p}{t}}(v_{\www})}^{\frac{1}{t}}
\\ 
&\lesssim  \left\|\mathcal{M}^{\sharp}\left(\big|T_{[a]}(f_1,\dots,f_l)\big|^t\right)
\right\|_{L^{\frac{p}{t}}(v_{\www})}^{\frac{1}{t}}
= \left\|\left(\mathcal{M}^{\sharp}\left(\big|T_{[a]}(f_1,\dots,f_l)\big|^t\right)\right)^{\frac{1}{t}}
\right\|_{L^{p}(v_{\www})},
\end{align*}
where $\mathcal{M}^{\sharp}=\mathcal{M}^{\sharp}_1$ and in the second inequality
we used the fact that $\|\mathcal{M}f\|_{L^{p_0}(w)} \lesssim \|\mathcal{M}^{\sharp}f\|_{L^{p_0}(w)}$
if $\mathcal{M}f \in L^{p_0}(w)$ and $w \in A_{\infty}$ (see \cite[Chapter IV, Theorem 2.20]{Ga_Ru1985}).
We observe that
$\big( \mathcal{M}^{\sharp}\big(|f|^t\big)(x)\big)^{\frac{1}{t}}
\le \mathcal{M}^{\sharp}_{t} f(x)$ for $0<t<1$,
since
$\left| |f(y)|^t - |c_Q|^t \right| \le |f(y) - c_Q|^t$.
Therefore, using H\"older's inequality with $t<\frac{r}{l}$ and our assumption \eqref{ass-smfe},
we have
\begin{align*}
\left( \mathcal{M}^{\sharp}\left(\big|T_{[a]}(f_1,\dots,f_l)\big|^t\right)(x)\right)^{\frac{1}{t}}
&\le \mathcal{M}^{\sharp}_{t}\big(T_{[a]}(f_1,\dots,f_l)\big)(x)
\\
&\le \mathcal{M}^{\sharp}_{\frac{r}{l}}\big(T_{[a]}(f_1,\dots,f_l)\big)(x)
\lesssim \mathbf{M}_q\big(f_1,\dots,f_l\big)(x),
\end{align*}
which together with \eqref{Le_Om_Pe_To_Tr2009_result} implies
\[
\big\|T_{[a]}(f_1,\dots,f_l) \big\|_{L^p(v_{\www})}
\lesssim \big\|\mathbf{M}_q(f_1,\dots,f_l)\big\|_{L^p(v_{\www})}
\lesssim_{a,\www}\prod_{j=1}^l \|f_j\|_{L^{p_j}(w_j)}.
\]
The proof is complete.

\hfill

%\bibliographystyle{apalike}
%\bibliography{HeBib}

%\printindex

\end{document}